\documentclass[a4paper,11pt]{article}

\usepackage{latexsym}
\usepackage{amssymb}        
\newtheorem{thm}{Theorem}[section]
\newtheorem{lem}[thm]{Lemma}

\newtheorem{prop}[thm]{Proposition}

\input epsf

\title{Paperfolding and Catalan numbers}
\author{Roland Bacher
}
\begin{document}
\maketitle

{\it Abstract\footnote{Math. Class: 11B65, 11B85, 11C20, 15A23.
}:} This paper
reproves a few results concerning 
paperfolding sequences using properties of 
Catalan numbers modulo 2. The guiding principle can be described by:
Paperfolding $=$ Catalan modulo 2 $+$ signs given by $2-$automatic 
sequences.

\section{Main results}

In this paper, a {\it continued fraction} is a formal expression
$$b_0+\frac{a_1}{b_1+\frac{a_2}{b_2+\dots}}=b_0+a_1\Big/\Big(b_1+
a_2\Big/\Big(b_2+\dots\Big)\Big)$$
with $b_0,b_1,\dots,a_1,a_2,\dots\in\mathbb{C}[x]$ two sequences of 
complex rational functions. We denote such a continued fraction by
$$b_0+\frac{a_1\vert}{\vert b_1}+\frac{a_2\vert}{\vert b_2}+\dots
+\frac{a_k\vert}{\vert b_k}+\dots$$
and call it {\it convergent} if the coefficients of the (formal) Laurent-series
expansions of its {\it partial convergents}
$$\frac{P_0}{Q_0}=\frac{b_0}{1},\ \frac{P_1}{Q_1}=\frac{b_0b_1+a_1}{b_1},\dots,
\frac{P_k}{Q_k}=b_0+\frac{a_1\vert}{\vert b_1}+\frac{a_2\vert}{\vert b_2}
+\dots+\frac{a_k\vert}{\vert b_k}\ ,\dots$$
are ultimately constant.
The limit $\sum_{k\geq n} l_kx^k\in\mathbb{C}[[x]][x^{-1}]$
is then the obvious formal Laurent power series whose first coefficients
agree with the Laurent series expansion of $P_n/Q_n$ for $n$ huge enough.

Let $W_1$ be the word $(-x_1)\ x_1$ of length $2$ in the alphabet
$\{x_1,-x_1\}$. For $k\geq 2$, we define recursively a word $W_k$ of
length $2(2^k-1)$ over the alphabet $\{\pm x_1,\pm x_2,\dots,\pm x_k\}$
by considering the word
$$W_k=W_{k-1}\ x_k\ (-x_k)\ \overline{W_{k-1}}$$
constructed by concatenation where $\overline{W_{k-1}}$ denotes the reverse
of $W_{k-1}$ obtained by reading $W_{k-1}$ backwards. The first 
instances of the words $W_k$ are
$$\begin{array}{lcl}
W_1=(-x_1)\ x_1\\
W_2=(-x_1)\ x_1\ x_2\ (-x_2)\ x_1\ (-x_1)\\
W_3=(-x_1)\ x_1\ x_2\ (-x_2)\ x_1\ (-x_1)\ x_3\ (-x_3)\
(-x_1)\ x_1\ (-x_2)\ x_2\ x_1\ (-x_1)\end{array}$$
and these words are initial subwords of a unique, well-defined
infinite word $W_\infty=w_1\ w_2\ w_3,\dots$ 
with letters $w_i\in \{\pm x_1,\pm x_2,\dots\}$.

The following result has been known for some time, cf \cite{K} or 
\cite {Sh2}. Similar or equivalent
equivalent statements are for instance given in
\cite{MPSa}, Theorem 6.5.6 in \cite{AS}. 
It is closely related to paperfolding, see \cite{AS}, \cite{PS}, \cite{MD}.
Related results are contained in \cite{LS}, \cite{SW} and \cite{Sh2}.

\begin{thm}\label{main} Consider a sequence $x_1,x_2,\dots $ with values
in ${\mathbb C}$. If the corresponding continued fraction
$$\frac{1\vert}{\vert 1}+\frac{w_1\vert}{\vert 1}+\frac{w_2\vert}{\vert 1}+
\dots+\frac{w_k\vert}{\vert 1}+\dots$$
associated with the sequence $W_\infty =w_1\ w_2\ \dots$
defined by the letters of the infinite word $W_\infty$
converges, then its limit is given by
$$1+x_1+x_1^2x_2+x_1^4x_2^2x_3+x_1^8x_2^4x_3^2x_4+\dots=
1+\sum_{j=1}^\infty\prod_{i=1}^j x_i^{2^{j-i}}\ .$$
\end{thm}

{\bf Remarks.} (i) In fact, the proof will follow from the identity
$$\frac{P_{2(2^k-1)}}{Q_{2(2^k-1)}}=
\frac{1\vert}{\vert 1}+\frac{w_1\vert}{\vert 1}+\frac{w_2\vert}{\vert 1}+
\dots+\frac{w_{2(2^k-1)}\vert}{\vert 1}=1+
\sum_{j=1}^k\prod_{i=1}^j x_i^{2^{j-i}}$$
(cf. Lemma 6.5.5 in \cite{AS})
showing that the $2(2^k-1)-$th partial convergent of this
continued fraction is a polynomial in $x_1,\dots,x_k$.

\ \ (ii) It is of course possible to replace the field 
$\mathbb{C}$ by any other reasonable field or ring.

{\bf Examples.} (1) Setting $x_i=x$ for all $i$ and multiplying by
$x$ we get
$$\frac{x\vert}{\vert 1}-\frac{x\vert}{\vert 1}+\frac{x\vert}{\vert 1}
+\frac{x\vert}{\vert 1}-\frac{x\vert}{\vert 1}+\frac{x\vert}{\vert 1}
+\dots$$
$$=x+x^2+x^4+x^8+\dots=\sum_{k=0}^\infty x^{2^k}\ .$$
The sequence of signs $w_1,w_2,\dots\in\{\pm 1\}$ given by
$$\sum_{k=0}^\infty x^{2^k}=\frac{x\vert}{\vert 1}+w_1\frac{x\vert}{\vert 1}+
w_2\frac{x\vert}{\vert 1}+w_3\frac{x\vert}{\vert 1}+\dots$$
can be recursively defined by
$$\begin{array}{l}
\displaystyle
w_{4i+1}=-w_{4i+2}=(-1)^{i+1}\ ,\\
\displaystyle 
w_{8i+3}=-w_{8i+4}=(-1)^i\ ,\\
w_{8i+7}=-w_{8i+8}=w_{4i+3}\end{array}$$
and is $2-$automatic (cf. \cite{AS}).

\ \ (2) Setting $x_i=x^{1+3^{i-1}}$ and multiplying by $x$ we get
$$\frac{x\vert}{\vert 1}-\frac{x^2\vert}{\vert 1}+\frac{x^2\vert}{\vert 1}
+\frac{x^4\vert}{\vert 1}-\frac{x^4\vert}{\vert 1}+\frac{x^2\vert}{\vert 1}
+\dots$$
$$=x+x^3+x^9+x^{27}+\dots=\sum_{k=0}^\infty x^{3^k}\ .$$

\ \ (3) Setting $x_i=x^{1+(i-1)i!}$ and multiplying by $x$ we get
$$\frac{x\vert}{\vert 1}-\frac{x\vert}{\vert 1}+\frac{x\vert}{\vert 1}
+\frac{x^3\vert}{\vert 1}-\frac{x^3\vert}{\vert 1}+\frac{x\vert}{\vert 1}
+\dots$$
$$=x+x^2+x^6+x^{24}+\dots=\sum_{k=1}^\infty x^{k!}\ .$$

Consider the sequence $\mu=(\mu_0,\mu_1,\dots)=
(1,1,0,1,0,0,0,1,\dots)$ defined by
$x\sum_{i=0}^\infty \mu_ix^i= \sum_{k=0}^\infty x^{2^k}$. We associate 
to $\mu=(\mu_0,\mu_1,\dots)$ its {\it Hankel matrix}
$H$ whose $(i,j)-$th coefficient $h_{i,j}=\mu_{i+j},\ 0\leq i,j$ depends 
only on the sum $i+j$ of its indices and is given by $(i+j)-$th element of
the sequence $\mu$.

Consider the ($2-$automatic) sequence  
$s=(s_0,s_1,s_2,\dots)=1\ 1\ (-1)\ 1\ 
(-1)\ (-1)\dots$ defined recursively by $s_0=1$ and
$$\begin{array}{l}
\displaystyle s_{2i}=(-1)^i\ s_i\ ,\\
\displaystyle s_{2i+1}=s_i\ .\end{array}$$
It is easy to check that $s_n=(-1)^{B_0(n)}$
where $B_0(n)$ counts the number of bounded blocks of consecutive
$0'$s of the binary integer $n$ (or, equivalently, $B_0(n)$ 
equals the number of blocks $10$ appearing in the binary expansion
of the $n$).
E.g, $720=2^9+2^7+2^6+2^4$ corresponds to the binary
integer $1011010000$. Hence $B_0(720)=3$.
This description shows that we have $s_{2^{n+1}+n}=-s_n$ for $n<2^n$.
Since we have $B_0({2^{n+1}+n})=1+B_0(n)$ for $n<2^n$ and 
$B_0({2^{n+1}+n})=B_0(n)$ for $2^n\leq n<2^{n+1}$,
the sequence $s_0,s_1,\dots$ can also be constructed by
iterating the application
$$W_\alpha W_\omega\longmapsto W_\alpha W_\omega (-W_\alpha) W_\omega$$
where where $W_\alpha,W_\omega$ are the first and the second half of the
finite word $s_0 \dots s_{2^n-1}$.

Denote by $D_s$ the infinite diagonal matrix with diagonal entries $s_0,s_1,s_2,\dots$
and by $D_a$ the infinite diagonal matrix with diagonal entries given by the 
alternating, $2-$periodic sequence $1,-1,1,-1,\overline{1,-1}$.

We introduce the infinite
unipotent lower triangular matrix $L$ with coefficients
$l_{i,j}\in \{0,1\}$ for $0\leq i,j$ given by 
$$l_{i,j}\equiv {2i\choose i-j}-{2i\choose i-j-1}
\equiv {2i\choose i-j}+{2i\choose i-j-1}\equiv 
{2i+1\choose i-j}\pmod 2$$
and denote by $L(k)$ the $k\times k$ submatrix with
coefficients $l_{i,j},\ 0\leq i,j<k$ of $L$. The matrix $L(8)$ for instance 
is given by
$$L(8)=\left(\begin{array}{rrrrrrrrrrrrrr}
1\\
1&1\\
0&1&1\\
1&1&1&1\\
0&0&0&1&1\\
0&0&1&1&1&1\\
0&1&1&0&0&1&1\\
1&1&1&1&1&1&1&1\\
\end{array}\right)\ .$$

\begin{thm} \label{LUdecomp} (i) We have
$$H=D_s\ L\ D_a\ L^t\ D_s\ .$$

\ \ (ii) ($BA0BAB-$construction for $L$)
We have for $n\geq 0$
$$L(2^{n+2})=
\left(\begin{array}{cccc}A\\B&A\\0&B&A\\B&A&B&A\end{array}\right)$$
where $A$ and $B$ are the obvious submatrices of
size $2^n\times 2^n$ in $L(2^{n+1})$.
\end{thm}

{\bf Remarks.} (i) Assertion (ii) of Theorem \ref{LUdecomp}
yields an iterative construction for $L$ by starting with 
$L(2)=\left(\begin{array}{cc}1\\ 1&1\end{array}\right)$.

\ \ (ii) Conjugating $H$ with an appropriate diagonal $\pm 1-$matrix, we can 
get the the LU-decomposition of any Hankel matrix associated with
a sequence having ordinary generating function
$x^{-1}\left(\sum_{k=0}^\infty \epsilon_k x^{2^k}\right),\ 
\epsilon_0=1,\epsilon_1,\epsilon_2,\dots\in \{\pm 1\}$. More precisely,
the appropriate conjugating diagonal matrix has coefficients
$d_{n,n}=(-1)^{\sum_{j=0} \nu_j\epsilon_{j+1}}$ for 
$0\leq n=\sum_{j=0} \nu_j2^j$ a binary integer.
The case $\epsilon_0=-1$ can be handled similarly by replacing 
the diagonal matrix $D_a$ with its opposite $-D_a$.
  
The matrix $P=\left(D_s\ L\ D_s\right)^{-1}=D_s\ L^{-1}\ D_s$ 
with coefficients $p_{i,j},\ 0\leq i,j$ encodes the formal 
orthogonal polynomials $Q_0,Q_1,\dots,Q_n=\sum_{k=0}^n
p_{n,k}x^k$ with moments 
$\mu_0,\mu_1,\dots$ given by $\sum_{i=0}^\infty \mu_ix^{i+1}=\sum_{k=0}^\infty
x^{2^k}$. 

In order to describe the matrix $P$ we introduce a unipotent 
lower triangular matrix $M$ with coefficients $m_{i,j}\in \{0,1\},\ 
0\leq i,j$ by setting
$$m_{i,j}={i+j\choose 2j}\pmod 2\ .$$ 
As above (with the matrix $L$) we denote by 
$M(k)$ the $k\times k$ submatrix with
coefficients $m_{i,j},\ 0\leq i,j<k$ of $M$.

\begin{thm} \label{Linverse} (i) We have 
$$L\ D_a\ M=D_a$$
(where $D_a$ denotes the diagonal matrix with $2-$periodic
diagonal coefficients $1,-1,1,-1,\dots$).

\ \ (ii) ($BA0BAB-$construction for $M$)
We have have for $n\geq 0$ 
$$M(2^{n+2})=
\left(\begin{array}{cccc}A'\\B'&A'\\A'&B'&A'\\B'&0&B'&A'\end{array}\right)$$
where $A'$ and $B'$ are the obvious submatrices of
size $2^n\times 2^n$ in $M(2^{n+1})$.
\end{thm}

The matrix $P=D_s\ D_a\ M\ D_a\ D_s=(D_s\ L\ D_s)^{-1}$ 
encoding the coefficients of 
the orthogonal polynomials with moments $\mu_0,\mu_1,\dots$ has 
thus all its coefficients in $\{0,\pm 1\}$, see also \cite{ALM} and
\cite{MS2}.

The (formal) orthogonal polynomials $Q_0,Q_1,\dots$ satisfy a classical 
three-term recursion relation with coefficients given by the continued 
Jacobi fraction expansion of $\sum \mu_i\ x^i=\frac{1}{x}\left(
\sum_{k=0}^\infty
x^{2^k}\right)$. This expansion is described by the following result.

\begin{thm} \label{Jacobifraction} We have
$$\begin{array}{rcl}
\displaystyle \sum_{k=0}^\infty x^{2^k}&=&\frac{x\vert}{\vert1-d_1x}+
\frac{x^2\vert}{\vert 1-d_2x}+\frac{x^2\vert}{\vert 1-d_3x}+\dots\\
\displaystyle &=&\frac{x\vert}{\vert1-x}+
\frac{x^2\vert}{\vert 1+2x}+\frac{x^2\vert}{\vert 1}
+\frac{x^2\vert}{\vert 1}+\frac{x^2\vert}{\vert 1-2x}+\dots\end{array}$$
where $d_n=\frac{s_n-s_{n-2}}{s_{n-1}}$ with $s_{-1}=0,\ s_0=1$
and $s_{2i+1}=(-1)^is_{2i}=s_i$ is the recursively defined
sequence introduced previously 
(extended by $s_{-1}=0$).
\end{thm}

The formal orthogonal polynomials with moments 
$\mu_0,\mu_1,\dots$ are thus 
recursively defined by $Q_0=1$,\ $Q_1=x+1$ and
$$Q_{n+1}=(x-d_n)Q_n+Q_{n-1}\hbox{ for }n\geq 2\ .$$

The rest of the paper is organised as follows.

Section 2 contains a proof of Theorem \ref{main}.

Section 3 recalls a few definitions and facts concerning Hankel
matrices, formal orthogonal polynomials etc.

Section 4 is a little digression about Catalan numbers.
All objects appearing in this paper are essentially obtained by
considering the corresponding objects associated with Catalan numbers
and by reducing them modulo $2$ using sign conventions
prescribed by recursively defined sequences (which are $2-$automatic).

Section 5 recalls two classical and useful results on
the reduction modulo $2$ (or modulo a prime number) of
binomial coefficients.

Section 6 is devoted to the proof of Theorem \ref{LUdecomp}.

Section 7 contains a proof of Theorem \ref{Linverse} and a related
amusing result.

Section 8 contains formulae for the matrix products
$M\ L$ and $L\ M$.

Section 9 contains a proof of Theorem \ref{Jacobifraction}.

Section 10 displays the $LU-$decomposition of the Hankel matrix
$\tilde H$ associated with the sequences $1,0,1,0,0,0,1,0,\dots$
obtained by shifting the sequence
$1,1,0,1,0,0,0,1,\dots$ associated with powers of $2$ by one (or equivalently,
by removing the first term).

Section 11 contains a uniqueness result.

\section{Proof of Theorem \ref{main}}

More or less equivalent reformulations of Theorem \ref{main} 
can be found in several places. A few are
\cite{MPSa}, \cite{PS} and \cite{AS}. 
We give her a computational proof: 
Symmetry-arguments of continued
fractions are replaced by polynomial identities.

Given a continued fraction
$$b_0+\frac{a_1\vert}{\vert b_1}+\frac{a_2\vert}{\vert b_2}+\dots$$
it is classical (and easy, cf. for instance
pages 4,5 of \cite{Pe}) to show that the product
$$\left(\begin{array}{cc}1&b_0\\0&1\end{array}\right)
\left(\begin{array}{cc}0&a_1\\1&b_1\end{array}\right)
\left(\begin{array}{cc}0&a_2\\1&b_2\end{array}\right)\cdots
\left(\begin{array}{cc}0&a_k\\1&b_k\end{array}\right)$$
equals
$$\left(\begin{array}{cc}P_{k-1}&P_k\\Q_{k-1}&Q_k\end{array}\right)$$
where 
$$\frac{P_k}{Q_k}=b_0+\frac{a_1\vert}{\vert b_1}+\dots+
\frac{a_k\vert}{\vert b_k}$$
denotes the $k-$th partial convergent.

Given a finite word $U=u_1\ u_2\ \dots\ u_l$ we consider the product
$$M(U)=\left(\begin{array}{cc}0&u_1\\1&1\end{array}\right)
\left(\begin{array}{cc}0&u_2\\1&1\end{array}\right)\cdots
\left(\begin{array}{cc}0&u_l\\1&1\end{array}\right)\ .$$
Defining
$$\begin{array}{lcl}
\displaystyle X_k&=&1+\sum_{j=1}^k\prod_{i=1}^j x_i^{2^{j-i}}\\
\displaystyle \tilde X_k&=&1+\sum_{j=1}^{k-1}\prod_{i=1}^j x_i^{2^{j-i}}-
\prod_{i=1}^k x_i^{2^{k-i}}\end{array}$$
(example: $X_3=1+x_1+x_1^2x_2+x_1^4x_2^2x_3$ and 
$\tilde X_3=1+x_1+x_1^2x_2-x_1^4x_2^2x_3$)
one checks easily the identities
$$X_n+x_{n+1}(X_{n-1}^2-X_n\tilde X_n)=X_{n+1}$$
and
$$X_n-x_{n+1}(X_{n-1}^2-X_n\tilde X_n)=\tilde X_{n+1}\ .$$

\begin{lem} \label{lemiracle} We have for all $n\geq 1$
$$M(W_n)=\left(\begin{array}{cc}
\tilde X_n-X_{n-1}^2&1-X_n\\
X_{n-1}^2& X_n\end{array}\right)$$
and
$$M(\overline{W_n})=\left(\begin{array}{cc}
X_n-X_{n-1}^2&1-\tilde X_n\\
X_{n-1}^2&\tilde X_n\end{array}\right)\ .$$
\end{lem}

\noindent{\bf Proof.} The computations
$$\begin{array}{cl}
\displaystyle
&M(W_1)=M((-x_1)\ x_1)=\left(\begin{array}{cc}0&-x_1\\1&1\end{array}\right)
\left(\begin{array}{cc}0&x_1\\1&1\end{array}\right)=
\left(\begin{array}{cc}-x_1&-x_1\\1&1+x_1\end{array}\right)\\
\displaystyle=
&\left(\begin{array}{cc}(1-x_1)-1^2&1-(1+x_1)\\1^2&1+x_1\end{array}\right)
\end{array}$$
and
$$M(\overline{ W_1})=M(x_1\ (-x_1))=\left(\begin{array}{cc}0&x_1\\1&1\end{array}\right)
\left(\begin{array}{cc}0&-x_1\\1&1\end{array}\right)=
\left(\begin{array}{cc}(1+x_1)-1^2&1-(1-x_1)\\1^2&1-x_1\end{array}\right)$$
prove Lemma \ref{lemiracle} for $n=1$.

For $n\geq 1$ we have
$$\begin{array}{cl}\displaystyle
&M(W_{n+1})=M(W_n\ x_{n+1}\ (-x_{n+1})\ \overline{ W_n})\\
\displaystyle
=&\left(\begin{array}{cc}\tilde X_n-X_{n-1}^2&1-X_n\\X_{n-1}^2&X_n\end{array}\right)
\left(\begin{array}{cc}x_{n+1}&x_{n+1}\\1&1-x_{n+1}\end{array}\right)
\left(\begin{array}{cc}X_n-X_{n-1}^2&1-\tilde X_n\\X_{n-1}^2&
\tilde X_n\end{array}\right)\\
\displaystyle
=&\left(\begin{array}{cc}(X_n-x_{n+1}(X_{n-1}^2-X_n\tilde X_n))-X_n^2&1-(
X_n+x_{n+1}(X_{n-1}^2-X_n\tilde X_n))\\
X_n^2&X_n+x_{n+1}(X_{n-1}^2-X_n\tilde X_n)\end{array}\right)\\
\displaystyle
=&\left(\begin{array}{cc}\tilde X_{n+1}-X_n^2&1-X_{n+1}\\X_n^2&X_{n+1}\end{array}\right)
\end{array}$$
and
$$\begin{array}{cl}\displaystyle
&M(\overline{ W_{n+1}})=M(W_n\ (-x_{n+1})\ x_{n+1}\ \overline{W_n})\\
\displaystyle
=&\left(\begin{array}{cc}\tilde X_n-X_{n-1}^2&1-X_n\\X_{n-1}^2&X_n\end{array}\right)
\left(\begin{array}{cc}-x_{n+1}&-x_{n+1}\\1&1+x_{n+1}\end{array}\right)
\left(\begin{array}{cc}X_n-X_{n-1}^2&1-\tilde X_n\\X_{n-1}^2&
\tilde X_n\end{array}\right)\\
\displaystyle
=&\left(\begin{array}{cc}(X_n+x_{n+1}(X_{n-1}^2-X_n\tilde X_n))-X_n^2&1-(
X_n-x_{n+1}(X_{n-1}^2-X_n\tilde X_n))\\
X_n^2&X_n-x_{n+1}(X_{n-1}^2-X_n\tilde X_n)\end{array}\right)\\
\displaystyle
=&\left(\begin{array}{cc}X_{n+1}-X_n^2&1-\tilde X_{n+1}\\X_n^2&
\tilde X_{n+1}\end{array}\right)
\end{array}$$
which proves Lemma \ref{lemiracle} by induction on $n$.\hfill$\Box$

{\bf Proof of Theorem \ref{main}.} Lemma \ref{lemiracle} shows that
$$\begin{array}{lcl}\displaystyle M(1\ W_n)&=&\left(\begin{array}{cc}
0&1\\1&1\end{array}\right)
\left(\begin{array}{cc}
\tilde X_n-X_{n-1}^2&1-X_n\\X_{n-1}^2&X_n\end{array}\right)\\
\displaystyle &=&
\left(\begin{array}{cc}
X_{n-1}^2&X_n\\
\tilde X_n&1\end{array}\right)\end{array}\ .$$
The $2(2^n-1)-$th partial convergent of the continued fraction
$$\frac{1\vert}{\vert 1}+\frac{w_1\vert}{\vert 1}+\frac{w_2\vert}{\vert 1}+
\dots$$ 
defined by the word $W_\infty$ equals thus 
$$X_n=1+\sum_{j=1}^n\prod_{k=1}^j x_k^{2^{j-k}}$$
and the sequence of these partial convergents has limit
$$X_\infty=1+\sum_{j=1}^\infty\prod_{k=1}^j x_k^{2^{j-k}}$$
which is the statement of Theorem \ref{main}.\hfill$\Box$

\section{Formal orthogonal polynomials}

The content of this section is classical and can for instance be found in
\cite{F} or \cite{V}.

We consider a formal power series
$$g(x)=\sum_{k=0}^\infty \mu_kx^k\in {\mathbb C}[[x]]$$
as a linear form on the complex vector space ${\mathbb C}[x]^*$
of polynomials by setting
$$l_g(\sum_{i=0}^n\alpha_i x^i)=\sum_{i=0}^n\alpha_i\mu_i\ .$$
We suppose that the $n-$th {\it Hankel matrix} $H(n)$ 
having coefficients
$h_{i,j}=\mu_{i+j},\ 0\leq i,j<n$ associated
with $\mu=(\mu_0,\dots)$ is invertible for all $n$. Replacing such a
sequence $\mu=(\mu_0,\mu_1,\dots)$ by
$\frac{1}{\mu_0}\mu=(1,\frac{\mu_1}{\mu_0},\dots)$ we may also
suppose that
$\mu_0=1$. We get a non-degenerate scalar product on 
${\mathbb C}[x]$ by setting
$$\langle P,Q\rangle=l_g(PQ)\ .$$
Applying the familiar Gram-Schmitt orthogonalisation we obtain a
sequence $Q_0=1,\ Q_1=x-a_0$,
$$Q_{n+1}=(x-a_n)Q_n-b_nQ_{n-1},\ n\geq 1$$
of mutually orthogonal monic polynomials, called {\it formal
orthogonal polynomials with moments $\mu_0,\mu_1,\dots$}.

We encode the above recursion relation for the polynomials
$Q_0,Q_1,\dots$ by the {\it Stieltjes matrix} 
$$S=\left(\begin{array}{cccccccccc}a_0&1\\b_1&a_1&1\\&b_2&a_2&1\\&&b_3&
a_3&1\\&&&\ddots&\ddots&\ddots\end{array}\right)\ .$$
 
The non-degeneracy condition on $\mu_0,\mu_1,\dots$
shows that the (infinite) Hankel matrix $H$ with coefficients 
$h_{i,j}=\mu_{i+j},\ 0\leq i,j$ has an $LU-$decomposition $H=LU$ where
$L$ is lower triangular unipotent and $U=DL^t$
is upper triangular with $D$ diagonal and invertible.

The diagonal entries of $D$ are given by
$$d_{i,i}=\prod_{k=1}^ib_k$$
and we have thus
$$\det(H(n))=\prod_{k=1}^{n-1}(b_k)^{n-k}$$
where $H(n)$ is the $n-$th Hankel matrix with coefficients
$h_{i,j}=\mu_{i+j},\ 0\leq i,j<n$.
The $i-$th row vector
$$(m_{i,0},m_{i,1},\dots)$$
of $M=L^{-1}$ encodes the coefficients of the $i-$th orthogonal 
polynomial
$$Q_i=\sum_{k=0}^i m_{i,k}x^k\ .$$
The Stieltjes matrix $S$ satisfies
$$LS=L_-=\left(\begin{array}{cccccc}
l_{1,0}&l_{1,1}\\
l_{2,0}&l_{2,1}&l_{2,2}\\
\vdots&&\ddots&\ddots\end{array}\right)$$
where $L_-$ is obtained by removing the first row of the 
unipotent lower triangular matrix $L$ with coefficients
$l_{i,j},\ 0\leq i,j$ involved in the $LU-$decomposition
of $H$.

The generating function $g=\sum_{k=0}^\infty \mu_k x^k$
can be expressed as a {\it continued fraction of Jacobi type}
$$g(x)=1\big/\big(1-a_0x-b_1x^2\big/\big(1-a_1x-b_2x^2\big/\big(
1-a_2x-b_3x^2\big/\big(1-\dots\big)\big)\big)\big)$$
where $a_0,a_1,\dots $ and $b_1,b_2,\dots$ are as above and are
encoded in the Stieltjes matrix.

\section{The Catalan numbers}

The Catalan numbers $C_n={2n\choose n}\frac{1}{n+1}$
are the coefficients of the algebraic generating function
$$\begin{array}{ccccl}
\displaystyle c(x)&=&\sum_{n=0}^\infty
C_nx^n&=&\frac{1-\sqrt{1-4x}}{2x}=1+x+2x^2+5x^3+14x^4+\dots\\
&&&=&\displaystyle\frac{1\vert}{\vert 1}-\frac{x\vert}{\vert 1}-
\frac{x\vert}{\vert 1}-\frac{x\vert}{\vert 1}-\dots\\
&&&=&\displaystyle\frac{1\vert}{\vert 1-x}-\frac{x^2\vert}{\vert 1-2x}-
\frac{x^2\vert}{\vert 1-2x}-\frac{x^2\vert}{\vert 1-2x}-\dots\end{array}$$
and satisfies $c=1+x\ c^2$.

We have $g(x)\equiv c(s)\pmod 2$ where $g(x)=\sum_{j=0}^\infty x^{2^j-1}$,
cf. the last lines of \cite{Sh}. A different proof can be given
by remarking that $c_n$ counts the number of planar binary rooted trees
with $n+1$ leaves. Call two such trees equivalent if they are
equivalent as abstract (non-planar) rooted trees. The cardinality of
each equivalence class is then a power of two. An equivalence
classe containing only one element is given by a rooted $2-$regular
trees having $2^n$ leaves which are all at the same distance $n$ 
from the root. 

Many interesting mathematical objects associated with
$g(x)$ are in fact obtained by reducing the corresponding objects of 
$c(x)$ modulo 2. Miraculously, all this works over $\mathbb Z$ by choosing
suitable signs given by a few recursively-defined sequences (which are in
fact always $2-$automatic, see \cite{AS} for a definition).

Consider the unipotent lower triangular matrices
$$L=\left(\begin{array}{rrrrrrr}
1\\
1&1\\
2&3&1\\
5&9&5&1\\
14&28&20&7&1\\
42&90&75&35&9&1
\end{array}\right)$$
and
$$\tilde L=\left(\begin{array}{rrrrrrr}
1\\
2&1\\
5&4&1\\
14&14&6&1\\
42&48&27&8&1\\
132&165&10&44&10&&1
\end{array}\right)$$
with coefficients 
$l_{i,j}={2i\choose i-j}-{2i\choose i-j-1}$, respectively
$\tilde l_{i,j}={2i+1\choose i-j}-{2i+1\choose i-j-1}$ for
$0\leq i,j$ obtained by taking all even (respectively odd)
lines (starting with line number
zero) of the so-called
{\it Catalan-triangle}
$$\begin{array}{rrrrrrrrrrrrr}
1\\
&1\\
1&&1\\
&2&&1\\
2&&3&&1\\
&5&&4&&1\\
5&&9&&5&&1\\
&14&&14&&6&&1\\
14&&28&&20&&7&&1
\end{array}$$
The matrix $L$ is associated with the $LU-$decomposition
of the Hankel matrix $H=L\ L^t$ having coefficients $h_{i,j}=
C_{i+j}={2(i+j)\choose i+j}/(i+j+1),\ 0\leq i,j$ the
Catalan numbers. The matrix $\tilde L$ yields the $LU-$decomposition
of the Hankel
matrix $\tilde H=\tilde L\ \tilde L^t$ having coefficients
$\tilde h_{i,j}=C_{i+j+1}={2(i+j+1)\choose i+j+1}/(i+j+2),\ 0\leq i,j$
associated with the shifted 
Catalan sequence $1,2,5,14,42,\dots$.

The inverse matrices $L^{-1}$ and $\tilde L^{-1}$ are given by $D_a\ M\ D_a$
and $D_a\ \tilde M\ D_a$ where $D_a$ is the diagonal matrix with 
2-periodic alternating diagonal coefficients $1,-1,1,-1,$ and where
$M$ and $\tilde M$ have
coefficients $m_{i,j}={i+j\choose 2j}$ and 
$\tilde m_{i,j}={i+j+1\choose 2j+1}$ for $0\leq i,j$.

Let us mention that the products $P=LM$ and $\tilde P=
\tilde L\ \tilde M$ have non-zero coefficients $p_{i,j}=
\frac{(2i)!\ j!}{i!\ (2j)!\ (i-j)!}$ and 
$\tilde p_{i,j}=4^{(i-j)}{i\choose j}$ for $0\leq j<i$. They 
can also be given by
$$P=L\ M=\hbox{exp}(\left(\begin{array}{ccccccccccc}
0\\2&0\\
&6&0&\\
&&10&0&\\
&&&\ddots&\ddots\end{array}\right))$$
and
$$\tilde P=\tilde L\ \tilde M=\hbox{exp}(\left(\begin{array}{ccccccccccc}
0\\4&0\\
&8&0&\\
&&12&0&\\
&&&\ddots&\ddots\end{array}\right))\ .$$

Computations suggest that the products $M\ L$ and $\tilde M\ \tilde L$
have also nice logarithms:
$$M\ L=\hbox{exp}(\left(\begin{array}{ccccccccccc}
0\\2&0\\
0&6&0&\\
2&0&10&0&\\
0&6&0&14&0\\
2&0&10&0&18&0\\
0&6&0&14&0&22&0\\
\vdots&\vdots&&&&&\ddots&\ddots\end{array}
\right))$$
and
$$\tilde M\ \tilde L=\hbox{exp}(\left(\begin{array}{ccccccccccc}
0\\4&0\\
0&8&0&\\
4&0&12&0&\\
0&8&0&16&0\\
4&0&12&0&20&0\\
\vdots&\vdots&&&&\ddots&\ddots
\end{array}\right))\ .$$

\section{Binomial coefficients modulo 2}

Chercher ref. (Concrete Maths? Graham, Knuth, Patashnik)

Using the Frobenius automorphism we have
$$(1+x)^{\sum_i \nu_i2^i}\equiv \prod_i (1+x^{2^i})^{\nu_i}\pmod 2$$
which shows the formula
$${n\choose k}\equiv \prod_i {\nu_i\choose \kappa_i}\pmod 2$$ 
where $n=\sum_{i\geq 0} \nu_i\ 2^i,\ \nu_i\in\{0,1\}$ and 
$k=\sum_{i\geq 0} \kappa_i\ 2^i,\ \kappa_i\in\{0,1\}$ are two natural 
binary integers. Kummer showed that the $2-$valuation $v_2$ of 
${n\choose k}$ (defined by $2^{v_2}\vert {n\choose k}$ and 
$2^{v_2+1}\not\vert {n\choose k}$) equals
the number of carry-overs when adding the two natural numbers $k$ and $n-k$,
written in base $2$ (cf. \cite{W}).

As an application, one (re)proves easily that $C_n={2n\choose n}\frac{1}{n+1}
\equiv 1\pmod 2$ if and only if $n+1\in\{2^k\}_{k\in{\mathbb N}}$ is
a power of two.
Indeed, ${2n\choose n}\frac{1}{n+1}\equiv \frac{(2n+1)!}{n!\ 'n+1)!}={2n+1\choose n}\pmod 2$
and addition of the binary integers $n$ and $(n+1)$ needs always a carry
except if $n=\sum_{i=0}^{k-1} 2^k$ and $n+1=2^k$.

\section{Proof of Theorem \ref{LUdecomp}}

The aim of this section is to prove Theorem \ref{LUdecomp}
which describes the $LU-$de\-composition 
of the Hankel matrix $H$ with coefficients $h_{i,j}=\mu_{i+j},\ 0\leq i,j$
where 
$$\sum_{k=0}^\infty \mu_ix^i=\sum_{j=0}^\infty x^{2^j-1}=1+x+x^3+x^7
+x^{15}+\dots.$$

We give in fact $3/2$ proofs: We first prove assertion (i) by 
using the definition $l_{i,j}\equiv {2i+1\choose i-j}\pmod 2$ of $L$.

We prove then assertion (ii) showing 
that the matrix $L$ satisfies the $BA0BAB-$construction.

Finally, we reprove assertion (i) using the recursive
$BA0BAB-$structure of $L$.

{\bf Proof of assertion (i) in Theorem \ref{LUdecomp}.} 
We denote by $\overline{n\choose k}\in \{0,1\}$
the reduction $\pmod 2$ with value in $\{0,1\}$ of the binomial coefficient
${n\choose k}$.
The main ingredient of the proof is the fact that
$$\overline{2n+1\choose k}=\overline{n\choose \lfloor k/2\rfloor}$$
where $\lfloor k/2\rfloor=k/2$ if $k$ is even and 
$\lfloor k/2\rfloor=(k-1)/2$ if $k$ is odd.

We compute first the coefficient $r_{2i,2j}$ of the product
$R=D_s L D_a L^t D_s$ (recall that $L$ has coefficients
$l_{i,j}=\overline{2i+1\choose i-j}$):
$$\begin{array}{lcl}
\displaystyle r_{2i,2j}&=&s_{2i}s_{2j}\sum_k (-1)^k\overline{4i+1\choose 2i-k}
\overline{4j+1\choose 2j-k}\\
\displaystyle &=&s_{2i}s_{2j}\left(\sum_{k=0}\overline{4i+1\choose 2i-2k}
\overline{4j+1\choose 2j-2k}-\sum_{k=1}
\overline{4i+1\choose 2i-2k+1}\overline{4j+1\choose 2j-2k+1}\right)\\
\displaystyle&=&s_{2i}s_{2j}\left(\sum_{k=0}\overline{2i\choose i-k}
\overline{2j\choose j-k}-\sum_{k=1}
\overline{2i\choose i-k}\overline{2j\choose j-k}\right)\\
&=&s_{2i}s_{2j}\overline{2i\choose i}\overline{2j\choose j}=\left\lbrace
\begin{array}{ll}
\displaystyle 1&\hbox{if }i=j=0\\
\displaystyle 0\qquad &\hbox{otherwise.}
\end{array}\right.
\end{array}$$

Similarly, we get
$$\begin{array}{lcl}
\displaystyle r_{2i+1,2j+1}&=&s_{2i+1}s_{2j+1}\sum_k (-1)^k\overline{4i+3\choose 2i+1-k}
\overline{4j+3\choose 2j+1-k}\\
\displaystyle &=&s_{2i+1}s_{2j+1}\left(\sum_{k=0}\overline{4i+3\choose 2i+1-2k}
\overline{4j+3\choose 2j+1-2k}-\sum_{k=0}
\overline{4i+3\choose 2i-2k}\overline{4j+3\choose 2j-2k}\right)\\
\displaystyle&=&s_{2i+1}s_{2j+1}\left(\sum_{k=0}\overline{2i+1\choose i-k}
\overline{2j+1\choose j-k}-\sum_{k=0}
\overline{2i+1\choose i-k}\overline{2j+1\choose j-k}\right)=0\end{array}$$

For $r_{2i,2j+1}=r_{2j+1,2i}$ we have
$$\begin{array}{lcl}
\displaystyle r_{2i,2j+1}&=&s_{2i}s_{2j+1}\left(
\sum_k \overline{4i+1\choose 2i-2k}\overline{4j+3\choose 2j+1-2k}
-\sum_k\overline{4i+1\choose 2i-2k-1}\overline{4j+3\choose 2j-2k}\right)
\\
&=&s_{2i}s_{2j+1}\left(\sum_k \overline{2i\choose i-k}\overline{2j+1
\choose j-k}
-\sum_k\overline{2i\choose i-k-1}\overline{2j+1\choose j-k}\right)\\
&=&s_{2i}s_{2j+1}\left(\sum_k \left(\overline{2i\choose i-k}+
\overline{2i\choose i-k-1}\right)\overline{2j+1\choose j-k}
-2\sum_k\overline{2i\choose i-k-1}\overline{2j+1\choose j-k}\right)\\
&=&s_{2i}s_{2j+1}\left(\sum_k \overline{2i+1\choose i-k}
\overline{2j+1\choose j-k}
-2\sum_k\overline{2i\choose i-k-1}\overline{2j+1\choose j-k}\right)\ .\\
\end{array}$$
Using parity arguments and the recursion 
$s_{2i}=(-1)^is_i,\ s_{2j+1}=s_j$ we get for
$i$ even
$$\begin{array}{lcl}
\displaystyle r_{2i,2j+1}&=&s_{2i}s_{2j+1}
\left(\sum_k(-1)^k\overline{2i+1\choose i-k}\overline{2j+1\choose j-k}+
\right.\\
\displaystyle&&\left.
\quad +2\sum_k\left(\overline{2i+1\choose i-2k-1}-\overline{2i\choose
i-2k-2}\right)\overline{2j+1\choose j-2k-1}\right)\\
\displaystyle&=&s_is_j\left(\sum_k(-1)^k\overline{2i+1\choose i-k}\overline{2j+1\choose j-k}+0\right)=r_{i,j}\ .\end{array}$$

Similarly, for $i$ odd we have 
$$\begin{array}{lcl}
\displaystyle r_{2i,2j+1}&=&-s_{2i}s_{2j+1}\left(\sum_k(-1)^k\overline{2i+1\choose i-k}\overline{2j+1\choose j-k}+\right.\\
\displaystyle&&\left.
\quad -2\sum_k\left(\overline{2i+1\choose i-2k}-\overline{2i\choose
i-2k-1}\right)\overline{2j+1\choose j-2k}\right)\\
\displaystyle&=&s_is_j\left(\sum_k(-1)^k\overline{2i+1\choose i-k}\overline{2j+1\choose j-k}+0\right)=r_{i,j}\ .\end{array}$$
Finally, the above computations together with
induction on $i+j$ show that $r_{i,j}=0$ except if $i+j=2^k-1$ for 
some $k\in\mathbb N$. In this case we get $r_{i,j}=1$ which proves the result.
\hfill $\Box$

{\bf Proof of assertion (ii) in Theorem \ref{LUdecomp}.}
Assume that $L(2^{n+1})$ is given by $\left(\begin{array}{cc}
A\\B&A\end{array}\right)$. (This holds clearly for $n=0$.)

We have for $0\leq i,j<2^{n+1}$
$${2(i+2^{n+1})+1\choose (i+2^{n+1})-(j+2^{n+1})}=
{2i+2^{n+2}+1\choose i-j}\equiv {2i+1\choose i-j}\pmod 2$$
since $2i+1<2^{n+2}$. This shows already that $L(2^{n+2})$ is of the form
$\left(\begin{array}{cc}\tilde A\\
\tilde B&\tilde A\end{array}\right)$
with $\tilde A=L(2^{n+1})$.

Consider now a binomial coefficient of the form
$${2(2^{n+1}+i)+1\choose (2^{n+1}+i)-j}={2^{n+2}+2i+1\choose 2^{n+1}+i-j}$$
with $0\leq i,j<2^{n+1}$ which is odd. This implies
$i+j+1\geq 2^{n+1}$ since otherwise $(2^{n+1}+i)-j<(2^{n+1}+i)+j+1<
2^{n+2}<2^{n+2}+2i+1$ and there must therefore be a carry when
adding the binary integers $((2^{n+1}+i)-j)$ and $((2^{n+1}+i)+j+1)$.

Two such binary integers $((2^{n+1}+i)-j)$ and $((2^{n+1}+i)+j+1)$ add thus
without carry if and only if the binary integers $(2^{n+2}+2^{n+1}+i-j)$
and $(2^{n+1}+i+j+1-2^{n+2})$ add without carry. This shows the
equality 
$${2(2^{n+1}+i)+1\choose (2^{n+1}+i)-j}\equiv {2(2^{n+1}+i)+1\choose (2^{n+1}
+i)-(2^{n+2}-1-j)}$$
for $0\leq i,j<2^{n+1}$.
Geometrically, this amounts to the fact that the block $\tilde B$ is
the vertical mirror of the block $\tilde A$ and this symmetry is 
preserved by the $BA0BAB-$construction.\hfill $\Box$

\noindent{\bf $BA0BAB-$proof of assertion (i).} We have
$$\left(\begin{array}{cccc}1\\1&1\\0&1&1\\1&1&1&1\end{array}\right)
\left(\begin{array}{cccc}1\\&-1\\&&1\\&&&-1\end{array}\right)
\left(\begin{array}{cccc}1&1&0&1\\&1&1&1\\&&1&1\\&&&1\end{array}\right)$$
$$=\left(\begin{array}{cccc}1&1&0&1\\1&0&-1&0\\0&-1&0&0\\1&0&0&0
\end{array}\right)$$
and conjugation by
$$\left(\begin{array}{cccc}1\\&1\\&&-1\\&&&1\end{array}\right)\hbox{ yields }
\left(\begin{array}{cccc}1&1&0&1\\1&0&1&0\\0&1&0&0\\1&0&0&0
\end{array}\right)\ .$$

The proof is now by induction using the $BA0BAB-$construction of $L$.
Writing $\tilde A=D_a\ A^t$ and $\tilde B=D_a\ B^t$ (with $D_a$ 
denoting the obvious submatrix of finite size of the previously infinite 
matrix $D_a$) we have
$$\left(\begin{array}{cccc}A\\B&A\\0&B&A\\B&A&B&A\end{array}\right)
\left(\begin{array}{cccc}\tilde A&\tilde B&0&\tilde B\\
&\tilde A&\tilde B&\tilde A\\
&&\tilde A&\tilde B\\
&&&\tilde A \end{array}\right)$$
$$=\left(\begin{array}{cccc}A\tilde A&A\tilde B&0&A\tilde B\\
B\tilde A&B\tilde B+A\tilde A&A\tilde B&B\tilde B+A\tilde A\\
0&B\tilde A&B\tilde B+A\tilde A&B\tilde A+A\tilde B\\
B\tilde A&B\tilde B+A\tilde A&A\tilde B+B\tilde A&2B\tilde B+2A\tilde A
\end{array}\right)\ .$$
We have $B\tilde B+A\tilde A=0$ by induction. 
The result is now correct up to signs if $A\tilde B=-B\tilde A$.

We conjugate now this last matrix by the diagonal matrix with
diagonal entries $s_0,\dots ,s_{2^{n+2}-1}$ which we can
write using blocks
of size $2^n$ in the form:
$$\left(\begin{array}{cccc} D_\alpha\\& D_\omega\\&& -D_\alpha \\
&&&D_\omega\end{array}\right)\ .$$

We get (after simplification using induction and the assumption 
$A\tilde B=-B\tilde A$)
$$\left(\begin{array}{cccc}
D_\alpha A\tilde AD_\alpha&D_\alpha A\tilde B D_\omega&0&D_\alpha A\tilde B
D_\omega\\
D_\omega B\tilde AD_\alpha&0&-D_\omega A\tilde B D_\alpha\\
0&-D_\alpha B\tilde AD_\omega\\
D_\omega B\tilde AD_\alpha\end{array}\right)\ .$$
Supposing the identity $A\tilde B=-B\tilde A$ we get
$$-D_\alpha B\tilde A D_\omega=D_\alpha A\tilde B D_\omega
\hbox{ and }
-D_\omega A\tilde B D_\alpha=
D_\omega B\tilde A D_\alpha$$
which have the correct form (one's on the antidiagonal and
zeros everywhere else) by induction.

We have yet to prove that $A\tilde B=-B\tilde A$.
Writing $A=\left(\begin{array}{cc}a\\b&a\end{array}\right)$,\ 
$B=\left(\begin{array}{cc}&b\\b&a\end{array}\right)$ and $\tilde A=D_a\ A^t,
\tilde B=D_a\ B^t$ we have 
$$A\tilde B=\left(\begin{array}{cc}a\\ b&a\end{array} \right)
\left(\begin{array}{cc}&\tilde b\\
\tilde b&\tilde a\end{array}\right)=\left(\begin{array}{cc} 0 &a\tilde b\\
a\tilde b&0\end{array}\right)$$
and
$$B\tilde A=\left(\begin{array}{cc}0&b\\ b&a\end{array} \right)
\left(\begin{array}{cc}\tilde a&\tilde b\\
0&\tilde a\end{array}\right)=\left(\begin{array}{cc} 0 &b\tilde a\\
b\tilde a&0\end{array}\right)$$
which proves $A\tilde B=-B\tilde A$ since $a\tilde b=-b\tilde a$
by induction.
\hfill $\Box$

\section{Proof of Theorem \ref{Linverse}}

Denote by $D_e$ the diagonal matrix with diagonal coefficients
$1,0,1,0,1,0,\dots,\overline{1,0,}$ and by $D_o=1-D_e$ the diagonal
matrix with coefficients $0,1,0,1,0,1,\dots,\overline{0,1,}$. 
We have thus $1=D_e+D_o$ and $D_a=D_e-D_o$. The main
ingredient for proving Theorem \ref{Linverse} is the following proposition
which is also of independent interest.

\begin{prop} \label{MDL} We have
$$M\ D_e\ L=A+D_e,\qquad 
M\ D_o\ L=A+D_o$$
where $A$ is the strictly lower triangular matrix with
coefficients $a_{i,j}=1$ if $i>j$ and $a_{i,j}=0$ otherwise.
\end{prop}

{\bf Proof.} We denote by $\alpha_{i,j}=(M\ D_e\ 
L)_{i,j}$ respectively $\beta_{i,j}=(M\ D_o\ 
L)_{i,j}$ the corresponding coefficients of both products.
The proof is by induction on $i+j$. Obviously, $\alpha_{i,j}=
\beta_{i,j}=0$ for $i<j$. 
The proof for $i=j$ is obvious since $M$ and $L$ are lower triangular
unipotent matrices.
Consider now for $i>j$
$$\alpha_{i,j}=\sum_k \overline{i+2k\choose 4k}
\overline{4k+1\choose 2k-j}\ .$$
(recall that $m_{i,k}=\overline{i+k\choose 2k}\equiv {i+k\choose 2k}\pmod 2$
and $l_{k,j}=\overline{2k+1\choose k-j}\equiv {2k+1\choose k-j}\pmod 2$).

Since 
$$\overline{2i+2k\choose 4k}=\overline{2i+1+
2k\choose 4k}\hbox{ and }\overline{4k\choose 2k-2j}=
\overline{4k+1\choose 2k-2j}=\overline{4k+1\choose 2k-(2j-1)}$$
we have 
$$\alpha_{i,j}=
\sum_k\overline {\lfloor i/2\rfloor+k\choose 2k}
\overline{2k\choose k-\lfloor(j+1)/2\rfloor}$$
$$=\sum_{k,\ k\equiv \lfloor(j+1)/2\rfloor\pmod 2} 
\overline{\lfloor i/2\rfloor+k\choose 2k}
\overline{2k+1\choose k-\lfloor(j+1)/2\rfloor}$$
$$=\left\lbrace
\begin{array}{ll}
\displaystyle \alpha_{\lfloor i/2\rfloor,\lfloor (j+1)/2\rfloor}
&\hbox{if }\lfloor (j+1)/2\rfloor\equiv 0\pmod 2\ ,\\
\displaystyle \beta_{\lfloor i/2\rfloor,\lfloor (j+1)/2\rfloor}
&\hbox{if }\lfloor (j+1)/2\rfloor\equiv 1\pmod 2\ .
\end{array}\right.$$
We have similarly
$$\beta_{i,j}=\sum_k\overline{i+2k+1\choose 4k+2}
\overline{4k+3\choose 2k+1-j}$$
$$=\sum
\overline{2\lfloor(i-1)/2\rfloor+1+2k+1\choose 4k+2}
\overline{4k+3\choose 2k+1-2\lfloor j/2\rfloor}$$
$$=\sum_{k,\ k\equiv \lfloor(i-1)/2\rfloor\pmod 2}
\overline{\lfloor(i-1)/2\rfloor+k+1\choose 2k+1}
\overline{2k+1\choose k-\lfloor j/2\rfloor}$$
$$=\sum_{k,\ k\equiv \lfloor(i-1)/2\rfloor\pmod 2}
\overline{\lfloor(i+1)/2\rfloor+k\choose 2k}
\overline{2k+1\choose k-\lfloor j/2\rfloor}$$
$$=\left\lbrace
\begin{array}{ll}
\displaystyle \alpha_{\lfloor (i+1)/2\rfloor,\lfloor j/2\rfloor}
&\hbox{if }\lfloor (i-1)/2\rfloor\equiv 0\pmod 2\ ,\\
\displaystyle \beta_{\lfloor (i+1)/2\rfloor,\lfloor j/2\rfloor}
&\hbox{if }\lfloor (i-1)/2\rfloor\equiv 1\pmod 2\ .
\end{array}\right.$$

except if $\lfloor(i-1)/2\rfloor=\lfloor j/2\rfloor=s$ where 
$$\sum_{k,\ k\equiv s\pmod 2}\overline{s+k+1\choose 2k+1}
\overline{2k+1\choose k-s}=1$$
since only $k=s$ yields a non-zero contribution.
This ends the proof.\hfill$\Box$

{\bf Proof of assertion (i) in Theorem \ref{Linverse}} Using the 
notations of Proposition \ref{MDL} we have 
$$M\ D_a\ L=(M\ D_e\ L)-(M\ D_o\ L)=(A+D_e)-(A+D_o)=D_a$$
which shows that $D_a\ M\ D_a=L^{-1}$ and ends the proof.\hfill $\Box$

{\bf Proof of assertion (ii) in Theorem \ref{Linverse}}
We have 
$$m_{i,j}=\overline{i+j\choose 2j}=\overline{i+j+2^{n+2}\choose
2j+2^{n+2}}=m_{i+2^{n+1},j+2^{n+1}}$$ for $i,j<2^{n+1}$.
This shows (together with induction)
that the lower right corner of $M(2^{n+2})$
is given by $\left(\begin{array}{cc}A'\\B'&A'\end{array}\right)$
and has thus the correct form.

In order to prove the recurrence formula for 
the lower left corner, we remark that $B'$
is by induction the horizontal mirror of $A'$. We have thus to show that 
$m_{i+2^{n+1},j}=m_{2^{n+1}-1-i,j}$ or equivalently that
$${2^{n+1}+i+j\choose 2j}\equiv {2^{n+1}-1-i+j\choose 2j}\pmod 2$$
for $0\leq i,j<2^{n+1}$. 

Consider
$$\frac{(2^{n+1}+i+j)(2^{n+1}+i+j-1)\cdots (2^{n+1}+i-j+1)}{(2j)!}$$
for $0\leq i,j<2^{n+1}$.
Since all terms of the numerator are $<2^{n+2}$ we have
$$\begin{array}{lcl}
\displaystyle 
{2^{n+1}+i+j\choose 2j}&\equiv&\frac{(2^{n+1}+i+j-2^{n+2})\cdots
(2^{n+1}+i-j+1-2^{n+2})}{(2j)!}\pmod 2\\
\displaystyle &\equiv&\frac{(2^{n+1}-i-j)(2^{n+1}-i-j+1)\cdots(2^{n+1}-i+j-1)}
{(2j)!}\pmod 2\\
\displaystyle &\equiv&{2^{n+1}-1-i+j\choose 2j}\pmod 2\end{array}$$
which ends the proof.\hfill$\Box$

{\bf $BA0BAB-$proof of assertion (i) in Theorem \ref{Linverse}}
Computing $L(4)\ D_a\ M(4)$ (with $D_a$ denoting a
diagonal matrix with alternating entries $1,-1,1,-1,\dots$ of
the correct size) we have
$$\left(\begin{array}{cccc}1\\1&1\\0&1&1\\1&1&1&1\end{array}\right)
\left(\begin{array}{cccc}1\\&-1\\&&1\\&&&-1\end{array}\right)
\left(\begin{array}{cccc}1\\
1&1\\1&1&1\\1&0&1&1\end{array}\right)=D_a\ .$$
Writing 
$$\left(\begin{array}{cccc}
\tilde A'\\\tilde B'&\tilde A'\\\tilde A'&\tilde B'&\tilde A'\\
\tilde B'&0&\tilde B'&\tilde A'\end{array}\right)=
D_a\ \left(\begin{array}{cccc}
A'\\B'&A'\\A'&B'&A'\\B'&0&B'&A'\end{array}\right)$$
we have
$$\left(\begin{array}{cccc}
A\\B&A\\0&B&A\\B&A&B&A\end{array}\right)\left(\begin{array}{cccc}
\tilde A'\\\tilde B'&\tilde A'\\\tilde A'&\tilde B'&\tilde A'\\
\tilde B'&0&\tilde B'&\tilde A'\end{array}\right)$$
$$=\left(\begin{array}{cccc}
A\tilde A'\\B\tilde A'+A\tilde B'&A\tilde A'\\
B\tilde B'+A\tilde A'&B\tilde A'+A\tilde B'&A\tilde A'\\
2B\tilde A'+2A\tilde B'&A\tilde A'+B\tilde B'&B\tilde A'+A\tilde B'&A\tilde A'
\end{array}\right)$$
$$=\left(\begin{array}{cccc}
D_a\\0&D_a\\
B\tilde B'+A\tilde A'&0&D_a\\
0&A\tilde A'+B\tilde B'&0&D_a
\end{array}\right)$$
by induction.

Since we have by induction $A\tilde A'=D_a$, it will be enough 
to show that $B\tilde B'=-D_a$. Writing 
$$B=\left(\begin{array}{cc}0&b\\b&a\end{array}\right)
\hbox{ and }\tilde B'=
\left(\begin{array}{cc}\tilde a'&\tilde b'\\\tilde b'&0\end{array}\right)$$
we get
$$B\tilde B'=\left(\begin{array}{cc}0&b\\b&a\end{array}\right)
\left(\begin{array}{cc}\tilde a'&\tilde b'\\\tilde b'&0\end{array}\right)=
\left(\begin{array}{cc}b\tilde b'&0\\b\tilde a'+a\tilde b'&b\tilde b'\end{array}\right)\ .$$
Since $b\tilde a'+a\tilde b'=0$ by the computation above and 
$b\tilde b'=-a\tilde a'$ by induction, we get the result.\hfill$\Box$

\section{$ML$ and $LM$}

As before, we denote by $L$ and $M$ the lower triangular unipotent
matrices with coefficients $l_{i,j},m_{i,j}\in\{0,1\}$ defined by
$$l_{i,j}\equiv {2i+1\choose i-j}\pmod 2,\ m_{i,j}\equiv{i+j\choose 2j}
\pmod 2$$
for $0\leq i,j$.

\begin{prop} \label{ML} The matrix $R=M\ L$ has coefficients
$$r_{i,j}=\left\lbrace\begin{array}{ll} 0\qquad &\hbox{if }i<j\\
1&\hbox{if }i=j\\
2&\hbox{if }i>j
\end{array}\right.$$
\end{prop}

\begin{prop} \label{LM} The product $LM$ yields a matrix which can recursively 
be constructed by iterating
$$\left(\begin{array}{cc}A\\B&A\end{array}\right)
\longmapsto \left(\begin{array}{cccc}A\\B&A\\
2A&B&A\\
2B&2A&B&A\end{array}\right)$$
starting with $A=1$ and $B=2$.
\end{prop}

Assertion (i) of Theorem \ref{Linverse} implies of course easily the identities
$$(M\ L)^{-1}=D_a\ M\ L\ D_a\hbox{ and }
(L\ M)^{-1}=D_a\ L\ M\ D_a\ .$$

{\bf Proof of Proposition \ref{ML}}. This follows immediately
from Proposition \ref{MDL} since
$ML=(M\ D_e\ L)+(M\ D_o\ L)$.\hfill $\Box$ 

{\bf Proof of Proposition \ref{LM}.} 
Using the $BA0BAB-$construction of $L$ and $M$
we have
$$\left(\begin{array}{cccc}A\\B&A\\0&B&A\\B&A&B&A\end{array}\right)
\left(\begin{array}{cccc}A'\\
B'&A'\\A'&B'&A'\\B'&0&B'&A'\end{array}\right)$$
$$=\left(\begin{array}{cccc}AA'\\BA'+AB'&AA'\\
BB'+AA'&BA'+AB'&AA'\\
2BA'+2AB'&AA'+BB'&BA'+AB'&AA'\end{array}\right)$$
Computing
$$AA'=\left(\begin{array}{cc}a\\b&a\end{array}\right)
\left(\begin{array}{cc}a'\\b'&a'\end{array}\right)=
\left(\begin{array}{cc}aa'\\
ba'+ab'&aa'\end{array}\right)$$
and 
$$BB'=\left(\begin{array}{cc}0&b\\b&a\end{array}\right)
\left(\begin{array}{cc}a'&b'\\b'&0\end{array}\right)=
\left(\begin{array}{cc}bb'\\
ba'+ab'&bb'\end{array}\right)$$
shows $AA'=BB'$ by induction. This finishes the proof by induction.\hfill
$\Box$

\section{Proof of Theorem \ref{Jacobifraction}}

We prove Theorem \ref{Jacobifraction} by showing that
the $n\times n$ submatrix formed by the first 
$n$ rows and columns of the Stieltjes matrix associated with
$\sum_{k=0}^\infty x^{2^k}$ has determinant $s_n$ for all $n$.

{\bf Proof of Theorem \ref{Jacobifraction}}
Let $H$ and $L$ be as in Theorem \ref{LUdecomp}.
We denote by $L_-$ the infinite matrix obtained by deleting the first
row of $L$. Denote by $C$ the matrix  with coefficients $c_{i,j}=1$ if
$(i-j)\in\{-1,0\}$ and $c_{i,j}=0$ otherwise (for $0\leq i,j$).

Elementary properties of the Catalan triangle show now that
$L=L_-\ C$. This implies that the finite matrix $L_-(n)$
formed by the first $n$ rows and columns of $L_-$ has
determinant $1$. Denoting by $D_{s_+1}(n)$ the $n\times n$ dia\-gonal matrix 
with entries $s_1,\dots,s_n$ we see that the finite 
Stieltjes matrix $S(n)=D_s(n)\ (L(n))^{-1}D_s(n)\ D_{s_+1}(n)\ 
\ L_-(n)\ D_{s}(n)$ associated with $H$ has determinant $s_n\in\pm 1$. 

Easy computations show that $\det(H(2))=-1$ and $\det(H(2^k))=1$
for $0\leq k\not=1$. We have also $\det(H(2^k+a))=(-1)^a\ \det(H(2^k-a))$
for $0\leq a<2^k$. This shows that $\det(H(n))=(-1)^{n\choose 2}$ 
and implies that the Stieltjes matrix associated with
$H$ is of the form
$$S(n)=\left(\begin{array}{cccccccccccc}
1&1\\
-1&d_2&1\\
&-1&d_3&1\\
&&-1&d_4&1\\
&&&\ddots&\ddots&\ddots\\
\end{array}\right)\ .$$
We have thus
$$s_n=\det(S(n))=d_n\det(S(n-1))+\det(S(n-2))=d_ns_{n-1}+s_{n-2}$$
for $n>1$ which shows
$$d_n=\frac{s_n-s_{n-2}}{s_{n-1}}\in\{0,\pm 2\},\ n>1$$
and proves the result.\hfill$\Box$

\section{The shifted Hankel matrix}

We consider the shifted Hankel matrix $\tilde H$ with coefficients
$\tilde h_{i,j}\in\{0,1\},\ 0\leq i,j$ given by
$\tilde h_{i,j}=1$ if $i+j+2=2^k$ for some natural integer $k\geq 1$
and $\tilde h_{i,j}=0$ otherwise.
We describe the $LU-$decomposition of $\tilde H$. The
associated Stieltjes matrix can be recovered from Theorem \ref{main}
by setting $x_1=x_2=x_3=\dots=x^2$.

Define lower triangular matrices $\tilde L$ and $\tilde M$
with coefficients $\tilde l_{i,j},\ \tilde m_{i,j}\in \{0,1\}$ given by
$\tilde l_{i,j}\equiv {2i+1\choose i-j}-{2i+1\choose i-j-1}
\equiv {2i+2\choose i-j}\pmod 2$ and 
$\tilde m_{i,j}\equiv {i+j+1\choose 2j+1}\pmod 2$.

It easy to see that one has recursive formulae
$$\tilde l_{2i,2j}=l_{i,j}\equiv {2i+1\choose i-j}\pmod 2,\quad 
\tilde l_{2i+1,2j+1}=\tilde l_{i,j}$$ 
and 
$$\tilde m_{2i,2j}=m_{i,j}\equiv {i+j\choose 2j}\pmod 2,\quad 
\tilde m_{2i+1,2j+1}=\tilde m_{i,j}\ .$$

The products $\tilde L\ \tilde M$ and $\tilde M\ \tilde L$
satisfy analogous recursive formulae. 

In order to state the main result of this section we need also
two sign-sequences $\tilde s_0,\tilde s_1,\dots$ and $\tilde t_0,\tilde t_1$
with $\tilde s_i,\tilde t_i\in\{\pm 1\}$. The sequence $\tilde s_i$ is
recursively defined by
$$\tilde s_{2i}=(-1)^i\hbox{ and }\tilde s_{2i+1}=\tilde s_i$$
and is obtained by removing the initial term of
the famous paperfolding sequence
$$1,1,-1,1,1,-1,-1,1,1,1,-1,-1,1,-1,-1,\dots$$
describing the peaks and valleys in a strip of paper which has iteratedly 
been folded (with all foldings executed in a similar way), see e.g.
\cite{ALM}, \cite{AS} or \cite{MS2}.

The sequence $\tilde t_i$ is defined by $\tilde t_0=1$ and
$$\tilde t_{2i+1}=\tilde t_i,\ \tilde t_{4i}=(-1)^i\tilde t_{2i},\ 
\tilde t_{4i+2}=\tilde t_{2i}\ .$$

As always, we denote by $D_{\tilde s}$, respectively $D_{\tilde t}$
the diagonal
matrices with diagonal entries $\tilde s_0,\tilde s_1,\dots $ respectively
$\tilde t_0,\tilde t_1,\dots$.

\begin{thm} \label{Hankelshift} (i) We have
$$D_{\tilde t}\ \tilde L\ D_{\tilde s}\ \tilde L^t\ D_{\tilde t}=\tilde H\ .$$

\ \ (ii) We have 
$$\tilde L\ D_{\tilde s}\ \tilde M=\tilde M\ D_{\tilde s}\ \tilde L=
D_{\tilde s}\ .$$
\end{thm}

{\bf Proof.} Let $r_{i,j}$ be the coefficient $(i,j)$ of the
matrix product $D_{\tilde t}\ \tilde L\ D_{\tilde s}\ \tilde L^t\ D_{\tilde t}
$. Since $\tilde l_{i,j}=0$ if $i\not\equiv j\pmod 2$ we have easily
that $r_{i,j}= 0$ if $i\not\equiv j\pmod 2$. For $r_{2i,2j}$ the computation 
reduces to assertion (i) of Theorem \ref{LUdecomp} since $\tilde l_{2i,2j}=
\overline{4i+2\choose 2i-2j}=\overline{2i+1\choose i-j}$ and $\tilde t_{2i}
=s_i$.

For $r_{2i+1,2j+1}$ we get
$$\begin{array}{cl}
\displaystyle&
\tilde s_{2i+1}\tilde s_{2j+1}\sum_k\overline{4i+4\choose 2i+1-2k-1}
\tilde t_{2k+1}\overline{4j+4\choose 2j+1-2k-1}\\
\displaystyle
=&\tilde s_i\tilde s_j\sum_k\overline{2i+2\choose i-k}\tilde t_k\overline{2j
+2\choose j-k}=r_{i,j}\end{array}$$
which proves assertion (i) by induction on $i+j$.

In order to prove assertion (ii), we denote by $r_{i,j}$ the
coefficient $(i,j)$ of the product
$\tilde L\ D_{\tilde s}\ \tilde M$. For $i\not\equiv j\pmod 2$ we
have $\tilde l_{i,j}=
\overline{2i+2\choose i-j}=\tilde m_{i,j}=\overline{i+j+1\choose 2j+1}=0$.
This shows that $r_{i,j}=0$ if $i\not\equiv j\pmod 2$.
For $r_{2i,2j}$ we have
$$\sum_k\overline {4i+2\choose 2i-2k}\tilde s_{2k}\overline{
2k+2j+1\choose 4j+1}$$
$$=\sum_k\overline {2i+1\choose i-k}(-1)^k\overline{
k+j\choose 2j}=\left\lbrace\begin{array}{ll}0&\hbox{if }i\not=j\\
(-1)^i&\hbox{if }i=j\end{array}\right.$$
by Theorem \ref{Linverse} which proves assertion (ii) in this case.

For $r_{2i+1,2j+1}$ we get
$$
\sum_k\overline{4i+4\choose 2i+1-2k-1}\tilde s_{2k+1}\overline{2k+2j+3\choose
4j+3}$$
$$=
\sum_k\overline{2i+2\choose i-k}\tilde s_{k}\overline{k+j+1\choose
2j+1}=r_{i,j}=\left\lbrace\begin{array}{ll}0&\hbox{if }i\not=j\\
\tilde s_i=\tilde s_{2i+1}&\hbox{if }i=j\end{array}\right.$$
which proves the result by induction on $i+j$.\hfill $\Box$

\subsection{$BA0BAB-$constructions for $\tilde L$ and $\tilde M$}

Denote by $\tilde L_0$ the infinite strictly lower triangular matrix
obtained by adding a first row of zeros to $\tilde L$.
We define also $\tilde M_0$ by adding a first row of zeros 
to $\tilde M$.

It can be shown that the matrix $\tilde L_0$ is obtained by iterating
$$\left(\begin{array}{cccc}
A\\
B&A\\
\end{array}\right)\longmapsto \left(\begin{array}{cccc}
A\\
B&A\\
0&B&A\\
B&A&B&A\end{array}\right)$$
where one starts with with $\left(\begin{array} {cc}A\\
B&A\end{array}\right)=\left(\begin{array} {cc}0&0\\
1&0\end{array}\right)$.

Similarly, $\tilde M_0$ is obtained by iterating
$$\left(\begin{array}{cccc}
A\\
B&A\\
\end{array}\right)\longmapsto \left(\begin{array}{cccc}
A\\
B&A\\
A&B&A\\
B&0&B&A\end{array}\right)$$
starting with $\left(\begin{array} {cc}A\\
B&A\end{array}\right)=\left(\begin{array} {cc}0&0\\
1&0\end{array}\right)$.

\section{Uniqueness}

Sequences with generating function
$$\sum_{k=0}^\infty \epsilon_k\ x^{2^k},\ \epsilon_i\in\{\pm 1\}$$
are bijectively related to the set of all paperfolding sequences
and can be characerized in terms of Hankel matrices as follows: 

\begin{prop} Let $s_0,s_1 ,\dots $ be a sequence with $s_i\in\{\pm 1,0\}$ such 
that $\det(H(n)),\det(\tilde H(n))=\pm 1$ for all $n$ where
$H(n)$ and $\tilde H(n)$ have coefficients
$h_{i,j}=s_{i+j}$ and $\tilde h_{i,j}=s_{i+j+1}$ for $0\leq i,j<n$.
Then $\sum_{k=0}^\infty s_k\ x^{k+1}=\sum_{k=0}^\infty \epsilon_k\ x^{2^k}$
for a suitable choice $\epsilon_0,\epsilon_1,\dots \in\{\pm 1\}$.
\end{prop}

{\bf Proof.} Computing $\det(H(n))\pmod 2$ determines $s_{2n-2}\pmod 2$ 
and computing $\det(\tilde H(n))\pmod 2$ determines $s_{2n-1}\pmod 2$.
This shows that 
$\sum_{k=0}^\infty s_k\ x^{k+1}\equiv \sum_{k=0}^\infty x^{2^k}\pmod 2$
and the signs can of course be chosen arbitrarily.
\hfill$\Box$

\bigskip
I thank J.P. Allouche, M. Mend\`es France and J. Shallit for helpful remarks and
comments.

Roland Bacher, Institut Fourier, UMR 5582,
Laboratoire de Math\'ematiques, BP 74, 38402 St. Martin d'H\`eres Cedex,
France, Roland.Bacher@ujf-grenoble.fr

\end{document}